\documentclass[12pt]{article}

\usepackage{latexsym,amssymb,amsmath}

\newtheorem{thm}{Theorem}
\newtheorem{cor}[thm]{Corollary}

\newtheorem{lemma}[thm]{Lemma}

\def\bb{\widetilde\lambda}

\def\eps{\varepsilon}
\def\e{\varepsilon}

\def\g{\gamma}

\def\il{\left<}
\def\ir{\right>}
\def\Lam{\Lambda}
\def\lam{\lambda}
\def\l{\lambda}
    
\def\N {{\mathbb N}}  

\def\RR {{\mathcal R}}
\newcommand{\reals} {{\mathbb R}}  
\def\ttt{\vartheta}
\def\tlam{\widehat{\lambda}}

\def\ee{\textrm{\tiny E}}      
\def\ww{\textrm{\tiny W}}

\title{Average Case Tractability of \\
Non-homogeneous Tensor Product Problems}
\author{ M. A. Lifshits, A. Papageorgiou, H. Wo\'zniakowski }

\begin{document}
\maketitle
\date

\begin{abstract}

We study $d$-variate approximation problems
in the average case setting with respect to a zero-mean Gaussian measure
$\nu_d$. Our interest is focused on measures 
having a structure of non-homogeneous 
linear  tensor product, where covariance kernel of $\nu_d$ is a product 
of univariate kernels, 
\[  
   K_d(s,t) = \prod_{k=1}^d \RR_k(s_k,t_k) \qquad \mbox{for}\ \ \ s,t\in [0,1]^d.
\]
We consider the normalized average error of algorithms 
that use finitely many evaluations of arbitrary linear functionals.
The information complexity is defined as the minimal number 
$n^{\rm avg}(\eps,d)$ of such evaluations  
for error in the $d$-variate case to be at most $\e$.
The growth of $n^{\rm avg}(\eps,d)$ as a function of $\eps^{-1}$ and $d$ 
depends on the eigenvalues of the covariance operator 
of   $\nu_d$ and determines whether a problem is tractable or not. 
Four types of tractability are studied and for each of them we 
find the necessary and 
sufficient conditions in terms of the eigenvalues of 
the integral operator with kernel $\RR_k$.

We illustrate our results by considering approximation problems
related to the product of Korobov kernels $\RR_k$. 
Each $\RR_k$ is characterized by a weight~$g_k$ and a
smoothness~$r_k$. We assume that weights are non-increasing
and smoothness parameters are non-decreasing. Furthermore they may be related,
for instance $g_k=g(r_k)$ for some non-increasing function $g$. 
In particular, we show that approximation 
problem is strongly polynomially tractable, i.e.,
$n^{\rm avg}(\eps,d)\le C\,\eps^{-p}$ for all $d\in \N, \eps\in (0,1]$,
where $C$ and $p$ are independent of $\e$ and $d$, iff
\[
  \liminf_{k\to\infty}\ \frac{\ln\,\frac 1{g_k}} {\ln k} >1.
\]
For other types of tractability we also show necessary and
sufficient conditions in terms of the sequences $g_k$ and $r_k$.
\end{abstract}

\section{Introduction}\label{s1}

Multivariate problems occur in many applications. 
They are defined on classes of functions of $d$ variables. 
Often the number of variables $d$ is large. 
Examples include problems in computational finance, statistics and physics.
These problems have been studied for different error criteria and 
in different settings including the worst and average case settings.
The cost of an algorithm solving a problem 
depends on the accuracy~$\e$ and the number of variables~$d$.
A problem is intractable if the cost of any algorithm 
is an exponential function
of~$\e^{-1}$ or~$d$.
Otherwise, the problem is tractable. Different types of 
tractable problems have been considered in the literature.
In fact, tractability of multivariate problems has been recently
a very active research area,  
see~\cite{NW08,NW10,NW12} and the references therein. 

More precisely, the information complexity $n(\eps,d)$ 
of a problem is the minimal number of information operations needed
by an algorithm to solve the problem  
with accuracy $\e$. The allowed information operations 
consist of function evaluations, or, more generally, 
of evaluations of arbitrary continuous linear functionals. 
We have
\begin{itemize}
\item \emph{weak} tractability if $n(\eps,d)$ is \emph{not}
exponential in $d$ and $\eps^{-1}$,
\item \emph{quasi-polynomial}
tractability if $n(\eps,d)$ is of order
$\exp(\,t\,(1+\ln\,d)(1+\ln\,\eps^{-1}))$,
\item \emph{polynomial} tractability if $n(\eps,d)$ is of order
$d^{\,q}\,\eps^{-p}$,
\item \emph{strong} polynomial tractability if $n(\eps,d)$ 
is of order $\eps^{-p}$.
\end{itemize}
The bounds above hold for all $d$ and all $\eps\in(0,1)$
with the parameters $t,q,p$ and the pre-factors
independent of $d$ and $\eps^{-1}$.

Strong polynomial tractability is the most challenging property.
Then the information complexity is bounded independently of $d$.
One may think that this property may hold only for trivial problems.
Luckily,  as we shall see, the opposite is sometimes true.

On the other hand, many multivariate problems are intractable. 
In particular, they suffer from the curse of dimensionality.
One way to vanquish the curse is to shrink the class of functions by 
introducing the weights that monitor the influence of successive
variables and groups of variables. For sufficiently fast decaying
weights not only we vanquish the curse but obtain strong
polynomial tractability; a survey of such results may be found again
in~\cite{NW08,NW10,NW12}.

The other way to vanquish the curse is by increasing the 
smoothness of functions with respect to the successive variables. 
This approach was taken recently in~\cite{PW10} for the worst case 
multivariate approximation in Korobov spaces.
In this paper we extend the approach of~\cite{PW10} to 
the average case setting and, in a much broader context, 
to tensor product Gaussian random fields. 
In this case we denote $n(\eps,d)=n^{\rm avg}(\eps,d)$
and restrict ourselves to information operations 
given by arbitrary continuous linear functionals
since the use of function values leads to the same results due 
to~\cite{HWW08} and Chapter 24 of~\cite{NW12}.

More precisely, we consider non-homogeneous linear
multivariate tensor product problems
in the average case with the normalized error criterion. 
The normalized error is used to measure the error of an algorithm
relative to the error of the zero algorithm.
A precise problem statement is given in Section \ref{s2}. 
The study of non-homogeneous case is necessary since
homogeneous linear multivariate tensor product problems
are intractable with this error criterion; see Chapter 6 in \cite{NW08}. 

In Section \ref{s3} we recall spectral
conditions for different types of tractability
in the average case and prove some new conditions.
The conditions are given in terms of the eigenvalues of the
covariance operator of the corresponding Gaussian measure.

In Section \ref{s4} these conditions are applied to non-homogeneous 
tensor product approximation problems. 
We equip the space of continuous real functions
defined on the $d$-dimensional unit cube $[0,1]^d$
with a zero-mean Gaussian measure with a covariance kernel of the form
\[  
   K_d(s,t) = \prod_{k=1}^d \RR_k(s_k,t_k), \qquad s,t\in [0,1]^d.
\]
Then $n^{\rm avg}(\eps,d)$ depends on spectral properties 
of the univariate integral operators with kernels $\RR_k$.
The main results of the paper , Theorems \ref{t:tp_62} -- \ref{t:tp_wt}, 
present spectral conditions for
polynomial, quasi-polynomial and weak tractability 
in this tensor product setting.

In Section \ref{s5} we illustrate these results for
Korobov kernels,  
$$
   \RR_k(x,y):=1+2\,g_k \,\sum_{j=1}^\infty
   j^{-2r_k}\,\cos(2\pi\,j(x-y)),  \qquad x,y\in[0,1],
$$ 
with varying smoothness parameters $r_k$ such that 
$$
   \tfrac12< r_1\le r_2 \le\cdots,
$$ 
and weight parameters $g_k$ such that
$$
   1 \ge g_1\ge g_2 \ge\cdots >0.
$$ 
The sequences $\{r_k\}$ and $\{g_k\}$ may be related. We may have
$$
g_k=g(r_k)
$$
for some non-increasing function $g:[\tfrac12,\infty)\to[0,1]$. 
The popular choice for Korobov space is to take $g(r)=(2\pi)^{-2r}$. 

It turns out that:
\begin{itemize}
\item Weak tractability holds iff 
$$
\lim_{k\to\infty}g_k=0.
$$
\item Quasi-polynomial tractability holds iff
$$
\sup_{d\in\N}\ \frac1{\max(1,\ln\,d)}\ \sum_{k=1}^dg_k
\,\max\left(1,\ln\,\frac1{g_k}\right)<\infty,
$$
under the assumption that $\liminf_{k\to\infty}r_k/\ln\, k>0$.

\item Polynomial tractability is equivalent to strong polynomial
  tractability.
\item Strong polynomial tractability holds iff
\[
   \rho_g:=\liminf_{k\to\infty}\ \frac{\ln\,\frac1{g_k}}{\ln k} >1.
\]
If this holds then $n^{\rm avg}(\eps,d)\le C\e^{-p}$ and the smallest $p$ is
$$
   \max\left(\frac2{2r_1-1},\frac2{\rho_g-1}\right).
$$
\end{itemize}
Other applications of our approach to 
tensor products problems are given in \cite{LPW1}
for covariance kernels corresponding to Euler and Wiener integrated processes.
We summarize the results of \cite{LPW1}
in Section \ref{s6} and compare them to those of the Korobov case that we 
study here.  By adjusting the weights $g_k$, the Korobov case 
behaves either like the Euler or Wiener case.

\section{Problem setting}\label{s2}

In this section we introduce multivariate problems in the average case
setting. We define the information complexity and 
the different notions of tractability.
More can be found in e.g.,~\cite{NW08} and~\cite{TWW88}.

Let $F_d$ be a Banach space of $d$-variate real functions defined on
a Lebesgue measurable set $D_d\subset\reals^d$. 
The space $F_d$ is equipped with a zero-mean Gaussian measure $\mu_d$
defined on Borel sets of $F_d$. We denote by
$C_{\mu_d}:F_d^{\,*}\to F_d$ the covariance operator of $\mu_d$.
Let $H_d$ be a Hilbert space with inner product and norm denoted
by $\il \cdot,\cdot\ir_{H_d}$ and $\|\cdot\|_{H_d}$, respectively.

We want to approximate a continuous linear operator
$$
S_d\,:\ F_d\,\to\,H_d.
$$
Let $\nu_d=\mu_d\,S_d^{-1}$ be the induced measure. 
Then $\nu_d$ is a zero-mean Gaussian measure on the Borel sets of $H_d$
with covariance operator $C_{\nu_d}:H_d\to H_d$ given by
$$
C_{\nu_d}=S_d\, C_{\mu_d} \, S^*_d\ 
$$
where $S^*_d: H_d\to F_d^*$ is the operator dual to $S_d$. 

Then $C_{\nu_d}$ is self-adjoint, nonnegative definite,
and has finite trace. Let 
$(\l_{d,j},\eta_{d,j})_{j=1,2,\dots}$ denote its eigenpairs 
$$
C_{\nu_d}\eta_{d,j}=\lambda_{d,j}\,\eta_{d,j} \ \ \ \ \
\mbox{with} \ \ \ \lambda_{d,1}\ge\lambda_{d,2}\ge\cdots.
$$
Then 
$$
{\rm trace}(C_{\nu_d})=\sum_{j=1}^\infty\lambda_{d,j}=
\int_{H_d}\|g\|^2_{H_d}\,\nu_d({\rm d}g)=
\int_{F_d}\|S_df\|^2_{H_d}\,\mu_d({\rm d}f)<\infty.
$$

We approximate $S_df$ for $f\in F_d$ by algorithms $A_n$ 
that use $n$ function evaluations or $n$ evaluations of arbitrary 
continuous linear functionals. It suffices to consider the case of arbitrary
continuous functionals since it is known that the results are roughly
the same for function values, see~\cite{HWW08} and Chapter~24 of~\cite{NW12}.  
Without essential loss of
generality, see e.g.,~\cite{NW08} as well as~\cite{TWW88}, 
we can restrict ourselves in the average case setting to
linear algorithms $A_n$ of the form
$$
A_n(f)=\sum_{j=1}^nL_j(f)\,g_j \ \ \ \ \  
\mbox{with} \ \ \ \ \ L_j\in F^{\,*}_d,\ \ g_j\in H_d.
$$
The average case error of $A_n$ is defined as
$$
e^{\rm avg}(A_n)=\bigg(\int_{F_d}\|S_df-A_n(f)\|_{H_d}^2\,
\mu_d({\rm   d}f)\bigg)^{1/2}.
$$

For a given $n$, it is well known that the algorithm $A_n$ that
minimizes the average case error is of the form
\begin{equation}\label{optalg}
A_n(f)=\sum_{k=1}^n\il S_df,\eta_{d,k}\ir_{H_d}\eta_{d,k},
\end{equation}
and its average case error is
\begin{equation}\label{avgerror}
e^{\rm avg}(A_n)=\bigg(\sum_{j=n+1}^\infty\lambda_{d,j}\bigg)^{1/2}.
\end{equation}
For $n=0$ we have the zero algorithm $A_0=0$. Its average case
error is called the initial error, and is given by the square-root of
the trace of the operator $C_{\nu_d}$, i.e., by~\eqref{avgerror} with
$n=0$.

The average case information complexity $n^{\rm avg}(\eps,d)$
is defined as the minimal~$n$ for which there is an algorithm 
whose average case
error reduces the initial error by a factor~$\eps$,
\begin{equation}\label{infocomp}
n^{\rm avg}(\eps,d)=\min\bigg\{\,n\ \bigg| \ \ 
\sum_{j=n+1}^\infty\lambda_{d,j}\le
\eps^2\,\sum_{j=1}^\infty\lambda_{d,j}\bigg\}.
\end{equation}

We present the definitions of four types of tractability
that will be studied in this paper.
Let $S=\{S_d\}_{d=1,2,\dots}$ denote a sequence of multivariate problems.  
We say that
\begin{itemize}
\item $S$ is \emph{weakly tractable} iff
$$
\lim_{\eps^{-1}+d\to\infty}\frac{\ln\,
\max\left(1,n^{\rm avg}(\eps,d)\right)}{\eps^{-1}+d}=0.
$$
\item $S$ is \emph{quasi-polynomially tractable} iff there are
positive numbers $C$ and $t$ such that
$$
n^{\rm avg}(\eps,d)\le C\,\exp\big(\,t\,(1+\ln\,d)\,(1+\ln\,\eps^{-1})\,\big)
\ \ \ \ \ \mbox{for all}\ \ \ d=1,2,\dots, \ \ \eps\in(0,1).
$$
The infimum of $t$ satisfying the bound above is called
the exponent of quasi-polynomial tractability and is denoted by 
$t^{\,\rm qpol-avg}$.
\item $S$ is \emph{polynomially tractable} iff there are 
non-negative numbers $C,q$ and $p$ such that
$$
n^{\rm avg}(\eps,d)\le C\,d^{\,q}\,\eps^{-p}
\ \ \ \ \ \mbox{for all}\ \ \ d=1,2,\dots, \ \ \eps\in(0,1).
$$ 
\item ${\rm APP}$ is \emph{strongly polynomially tractable} iff there are 
positive 
numbers $C$ and $p$ such that
$$
n^{\rm avg}(\eps,d)\le C\,\eps^{-p}
\ \ \ \ \ \mbox{for all}\ \ \ d=1,2,\dots, \ \ \eps\in(0,1).
$$ 
The infimum of $p$ satisfying the last bound is called the exponent
of strong polynomial tractability and is denoted by $p^{\rm \,str-avg}$.
\end{itemize}

Tractability can be fully characterized in terms of the 
eigenvalues $\lambda_{d,j}$. Necessary and sufficient
conditions on weak, quasi-polynomial,
polynomial and strong polynomial tractability 
can be found in Chapter~6 of~\cite{NW08} and 
Chapter~24 of~\cite{NW12}. In particular, $S$ is polynomially tractable 
iff there exist $q\ge0$ and $\tau\in(0,1)$ such that
\begin{equation}\label{poltract}
C:= \sup_{d\in\N}
\frac{\left(\sum_{j=1}^\infty\lambda_{d,j}^\tau\right)^{1/\tau}}
{\sum_{j=1}^\infty\lambda_{d,j}}\ d^{\, -q} <\infty.
\end{equation}
If so then
\begin{equation}\label{poltract2}
n^{\rm avg}(\eps,d)\le
C^{\tfrac{\tau}{1-\tau}}
\, d^{\tfrac{q\tau}{1-\tau}}\,\eps^{\tfrac{-2\tau}{1-\tau}}
\end{equation}
for all $d\in\N$ and $\eps\in(0,1)$. 
  
Furthermore, $S$ is strongly polynomially tractable
iff~\eqref{poltract} holds with $q=0$. The exponent of strong
polynomial tractability is
\begin{equation}\label{expostpol}
   p^{\,\rm str-avg}=\inf\left\{\frac{2\tau}{1-\tau}\ \bigg| \ \ \tau \
   \mbox{satisfies~\eqref{poltract} with $q=0$}\right\}.   
\end{equation}

\section{General Bounds} \label{s3}

We show bounds on $n^{\rm avg}(\eps,d)$
which we will use to derive necessary and sufficient conditions for
the four types of tractability. We first analyze an arbitrary problem
$\{S_d\}$ and then restrict our attention to non-homogeneous tensor
product problems.

We begin with a bound on $n^{\rm avg}(\eps,d)$ 
which from a probabilistic point of view
is an application of Chebyshev's inequality. 

\begin{lemma} For any $\eps\in(0,1),\,d\in\N$, $\tau\in (0,1)$ and
  $z>0$  
we have
\begin{equation}\label{cheb}
n^{\rm avg}(\eps,d)
\le 
\frac{\sum_{j=1}^\infty\lam_{d,j}^z}{\left(\sum_{j=1}^\infty\lam_{d,j}\right)^z}
\ 
\left[ \frac{\sum_{j=1}^\infty\lam_{d,j}^\tau}
{\left(\sum_{j=1}^\infty\lam_{d,j}\right)^\tau}\right]^{z/(1-\tau)}\,
\eps^{-2z/(1-\tau)}.
\end{equation}
\end{lemma}

{\bf Proof.} Let 
$b:=\left[\left(\sum_j \lam_{d,j} \right)\eps^{2}/
\left(\sum_j \lam_{d,j}^\tau\right)\right]^{1/(1-\tau)}$. Then
\[
   \sum_{ j: \lam_{d,j} \le b } \lam_{d,j} 
\le  \sum_{ j: \lam_{d,j}\le b} \lam_{d,j}^\tau b^{1-\tau}
   \le  \sum_{j} \lam_{d,j}^\tau b^{1-\tau} = \sum_{j} \lam_{d,j}\ \eps^2.
\]   
Hence,
\begin{eqnarray*}
n^{\rm avg}(\eps,d)&\le&
\#\{ j: \lam_{d,j}> b \} =
\#\{ j: \lam_{d,j}^z > b^z \} 
\le\frac{\sum_{j:\ \lam_{d,j}^z > b^z}\lam_{d,j}^z}{b^z}\\
&\le& \frac{\sum_j \lam_{d,j}^z}{b^z} 
   = \frac{\sum_j\lam_{d,j}^z}{\left(\sum_j\lam_{d,j}\right)^z}\ 
 \left[ \frac{\left(\sum_j \lam_{d,j}^\tau \right)}
{\left(\sum_j \lam_{d,j}\right)^\tau} \right]^{z/(1-\tau)}\
\eps^{-2z/(1-\tau)},
\end{eqnarray*}
as claimed. \hfill $\Box$
\bigskip

Note that~\eqref{cheb} immediately proves sufficiency of 
polynomial tractability conditions in~\eqref{poltract}.  
Furthermore, if we set $z=\tau$ then we obtain the estimate
\eqref{poltract2} with 
the exponent of strong polynomial tractability 
at most $2\tau/(1-\tau)$ 
for $\tau$ satisfying~\eqref{poltract} with $q=0$.

As we shall see now, the bound~\eqref{cheb} is also useful
when we consider quasi-polynomial tractability. 
In the rest of the paper we denote
$$
\ln_+d:=\max(1,\ln d).
$$

\begin{thm}\label{t:qpt} 
$S$ is quasi-polynomially tractable iff  there exists
$\delta\in(0,1)$ such that
\begin{equation}\label{cr_qpt}
   \sup_{d\in \N}\  
    \frac{\sum_{j=1}^\infty \lam_{d,j}^{1-\frac{\delta}{\ln_+d}} }
         {\left(\sum_{j=1}^\infty \lam_{d,j}\right)^{1-
\frac{\delta}{\ln_+d}}}
    <\infty. 
\end{equation}
\end{thm}

{\bf Proof.}  Sufficiency.
Apply~\eqref{cheb} with $\tau=1-\frac{\delta}{\ln_+d}\in(0,1)$ and
$z=1$. We obtain
\begin{eqnarray*}
n^{\rm avg}(\eps,d) &\le& 
\left[ \frac{\left(\sum_j \lam_{d,j}^{1-\frac{\delta}{\ln_+d}}\right)
\eps^{-2}}
{\left(\sum_j \lam_{d,j}\right)^{1-\frac{\delta}{\ln_+d}}} 
\right]^{\ln_+d/\delta}
\\
&\le& M_\delta^{\ln_+d/\delta} \eps^{-2\ln_+d/\delta}
=\exp\left(\frac{\ln M_\delta}{\delta} \ln_+d+ \frac 2\delta \ln_+d\  
\ln\,\eps^{-1} \right),
\end{eqnarray*}
where $M_\delta$ is the supremum in (\ref{cr_qpt}). We can 
rewrite the last estimate as
$$
n^{\rm avg}(\eps,d) \le
\exp\left(t(1+\ln\,d)\,(1+\ln\,\eps^{-1})\right)
$$
for $t=\delta^{-1}\,\max(2,\ln\,M_\delta)$.
This means that $S$ is 
quasi-polynomially tractable.

\medskip

Necessity. Assume now that $S$ is quasi-polynomially tractable, i.e., 
\[
n^{\rm avg}(\eps,d) \le C\,\exp\left(t(\ln d +1)(\ln\,\eps^{-1}+1)\right).
\]

We show that there is  $\delta\in(0,1)$ such that 
\[
\sup_{d\in \N}\   
\frac 
         {\sum_j \lam_{d,j}^{1-\frac{\delta}{\ln_+d}           }} 
{\left( \sum_j \lam_{d,j} \right)^{1-\frac{\delta}{\ln_+d}} }
<\infty. 
\]

Note that the last condition is invariant under multiplying 
the eigenvalues by a positive number, and so 
is the value of  $n^{\rm avg}(\eps,d)$. 
That is why we may assume that $\sum_j \lam_{d,j}=1$.
Quasi-polynomial tractability means that for all $\eps>0$ and $d\ge 1$
we have 
\[
\sum_{j\ge  C\,\exp\left(t(\ln d +1)(|\ln \eps|+1)\right)+1 }
   \lam_{d,j}\le \eps^2.
\]
Let $\eps:=e\, (n/C)^{\frac{-1}{t(\ln d+1)}}$. Then
\begin{equation}\label{e1}
     \sum_{j>n} \lam_{d,j}  \le   e^2 
(n/C)^{\frac{-2}{t(\ln d+1) }} := e^2 (n/C)^{-h}
\end{equation}
with $h=2/(t(1+\ln\,d)$.
\medskip

To avoid too small eigenvalues, we introduce a regularization  
\[
   \tlam_{d,j}:= \max\{ \lam_{d,j}, h\,j^{-1-h} \}.
\]
Note that~\eqref{e1} implies 
\begin{equation}\label{e2}
     \sum_{j> n} \tlam_{d,j}\le  \sum_{j> n}\lam_{d,j} +
\sum_{j> n} h j^{-1-h}  \le (e^2C^h+1) n^{-h} .
\end{equation}
Let
\[
  N_m=\{  j\in\N:\ 2^{m/h}\le j< 2^{(m+1)/h} \}, \qquad m=0,1,....
\]
Note that the structure of $N_m$ depends on $h$.
For any $\gamma\in (0,1)$ and any integer $m\ge 0$ we have
\begin{eqnarray*}
  \sum_{j\in N_m} \lam_{d,j}^{1-\gamma}  
  &\le& \sum_{j\in N_m} \tlam_{d,j}^{1-\gamma}
  \le \sum_{j\in N_m} \tlam_{d,j} \ 
[ \min_{j\in N_m} \tlam_{d,j}]  ^{-\gamma}
\\
  &\le& \sum_{j\ge 2^{m/h}} \tlam_{d,j} \
\left[  h \left(2^{(m+1)/h}\right)^{-1-h}  \right]  ^{-\gamma}
\\
(\textrm{by}\ (\ref{e2}) )
  &\le&  (e^2C^h+1) \  2^{-m}  \cdot      
h^{-\gamma}   2^{  \frac{\gamma(m+1)(1+h)}{h}} .
\end{eqnarray*}
For a fixed  $\delta\in (0,1)$, 
let $\gamma=\frac{\delta h}{1+h}$. We obtain
 \begin{eqnarray*}
  \sum_{j\in N_m} \lam_{d,j}^{1-\gamma}  
&\le&   (e^2C^h+1)   2^{-m}  \cdot      h^{-\gamma}   2^{ \delta(m+1)} 
\\
&\le&  (e^2C^h+1)2^{\delta}  \cdot  h^{-\gamma}   2^{ -(1-\delta) m} 
\le  (e^2C^h+1) 2^{\delta} 
\exp\left(|\ln h| \frac{\delta h}{h+1}\right) \   2^{ -(1-\delta) m}  .
\end{eqnarray*}
Since
\[
   \sup_{0<h\le \frac{2}{t}} |\ln h| h =: c(t) <\infty,
\]
it follows that
\[
 \sum_j \lam_{d,j}^{1-\gamma} 
=\sum_{m=0}^\infty \sum_{j\in N_m} \lam_{d,j}^{1-\gamma} 
\le 2 (e^2C^h+1) e^{c(t)} \sum_{m=0}^\infty  2^{ -(1-\delta) m} =: c(t,\delta).
\]
Note that 
\[
\gamma=\frac{\delta h}{1+h}= \frac{2\delta}{t(1+\ln d)+2}.
\]
Thus we have 
\[
     \sup_{d\in \N}\   
\sum_j \lam_{d,j}^{1-\frac{2\delta}{t(1+\ln d)+2}} <\infty.
\]
Let
$$
\delta^{\,\prime}:=\min_{d\in\N}\frac{2\delta\,\ln_+d}{t(1+\ln\,d)+2}
\le \frac{2\delta}{t+2}<1.
$$
Then 
$$
1-\frac{2\delta}{t(1+\ln\,d)+2}\le 1-\frac{\delta^{\,\prime}}{\ln_+d}
\ \ \ \ \ \mbox{for all}\ \ \ d\in\N,
$$
and 
\[
     \sup_{d\in \N}\,   
\sum_j \lam_{d,j}^{1-\frac{\delta^{\,\prime}}{\ln_+d}} <\infty,
\]
as required. This completes the proof.
\hfill $\Box$

\bigskip
Theorem~\ref{t:qpt} does not address the exponent $t^{{\rm qpol-avg}}$
of quasi-polynomial tractability. There is, however, 
the bound on the exponent presented in the first part of the proof, 
\begin{equation}\label{tbound}
   t^{{\rm qpol-avg}}\le \delta^{-1}\,\max(2,\ln\,M_\delta)
\end{equation}
for all $\delta\in(0,1)$ satisfying~\eqref{cr_qpt}. 

The presence of $M_\delta$ may seem artificial.
However, we now show that in general $M_\delta$ 
cannot be avoided in determining 
the exponent of quasi-polynomial tractability.  
Indeed, for $\delta\in(0,1)$, $M>1$ and $d\ge 1$  let 
$N=N(d,M,\delta):=\lfloor M^{\ln_+ d/\delta}\rfloor$
and consider the following eigenvalues 
\[
\lam_{d,j}:=
\begin{cases}  
               1& \ \ \ \mbox{for}\ \ \ j=1,2,\dots, N,
               \\
               0&\ \ \ \mbox{for}\ \ \   j>  N.
\end{cases}
\]
Then 
\[
M_\delta=
\sup_{d\in \N}\  
    \frac{\sum_{j=1}^\infty \lam_{d,j}^{1-\frac{\delta}{\ln_+d}} }
         {\left(\sum_{j=1}^\infty \lam_{d,j}\right)^{1-
\frac{\delta}{\ln_+d}}}
= \sup_{d\in \N}\  N(d,M,\delta)^{\delta/\ln_+d}
=\lim_{d\to\infty}\  N(d,M,\delta)^{\delta/\ln_+d}
=M,
\] 
Hence quasi-polynomial tractability holds and for any $\eps\in(0,1)$ we have
\[
   n^{{\rm avg}}(\eps,d)=\lceil(1-\eps^2)N\rceil 
   \le C \exp\left( t(1+\ln d)(1+\ln\,\eps^{-1})\right).
\]
It follows that
\[
  t\ge \lim_{\eps\to 1}
\lim_{d\to\infty} \frac{\ln \lceil(1-\eps^2)N(d,M,\delta)\rceil }{\ln d}
  =\frac {\ln M}{\delta}= \frac {\ln M_\delta}{\delta}.
\]
This justifies 
the presence of \ $\frac {\ln M_\delta}{\delta}$\ in the bound 
\eqref{tbound} for exponent of quasi-polynomial tractability.
However, we believe this bound is not always sharp.
\vskip 1pc

We now show that the necessary condition on quasi-polynomial
tractability can be simplified by eliminating the powers of
$1-\delta/\ln_+d$. The following lemma will be a convenient tool 
for establishing this fact.

\begin{lemma}\label{useful1}
Let $\Lambda_d=\sum_{j=1}^\infty\lambda_{d,j}$. 
For any $\gamma>0$ we have 
\begin{equation}\label{jen}
\Lam_d^{-1} \sum_{j=1}^\infty \lam_{d,j}^{1-\gamma} 
\ge
\exp\left(
- \gamma  \sum_{j=1}^\infty  \frac{\lam_{d,j} 
\, \ln\,\lam_{d,j}\   }{\Lam_d}\right).
\end{equation}
\end{lemma}

{\bf Proof}. \ Jensen's inequality states that 
for a convex function $\phi(\cdot)$ defined on a convex set $D$, 
non-negative weights $p_j$ satisfying $\sum_j p_j=1$, and 
any set of arguments $x_j$ from $D$ we have
\[  
\sum_j p_j \phi(x_j) \ge  \phi  \left( \sum_j p_j x_j\right). 
\]
We apply Jensen's inequality with
$p_j:= \frac{\lam_{d,j}}{\Lam_d}$, $x_j:= -\ln\,p_j$ and the function
$\phi(x)=e^{\gamma x}$ for $x\in D:=\reals$.  We obtain
\begin{eqnarray*}
&&\frac{\sum_j\lam_{d,j}^{1-\gamma}} {\Lam_d^{1-\gamma}}
= \sum_j p_j^{1-\gamma} = \sum_j p_j \exp(-\gamma\ln p_j)
=\sum_j p_j \phi(x_j) 
\\
&\ge&  \phi  \left( \sum_j p_j x_j\right)
= \exp  \left( \gamma \sum_j (- p_j \ln\, p_j )\right)
= \exp  \left( -\gamma \sum_j   p_j 
\left(\ln\,\lam_{d,j}-\ln\,\Lam_d\right) \right)
\\
&=& \Lam_d^\gamma \exp  \left( -\gamma \sum_j  p_j  \ln\,\lam_{d,j} \right)
= \Lam_d^\gamma \exp  
\left( -\gamma \sum_j  \frac{\lam_{d,j} \ln\,\lam_{d,j}}{\Lam_d} \right).
\end{eqnarray*}
This is equivalent to~\eqref{jen} and completes the proof. \hfill $\Box$

\medskip

We will see in the next section that the right-hand side
of~\eqref{jen} is convenient for tensor product problems.
We are ready to simplify the necessary conditions for
quasi-polynomial tractability. 

\begin{cor}\label{corqpol}
 If quasi-polynomial tractability holds then
\begin{equation}\label{nes_qpt}
\sup_{d\in \N} \  \frac{1} {\ln_+d} \ 
\sum_{j=1}^\infty  \frac{\lam_{d,j}}{\Lam_d} \, \ln\,
\left(\frac{\Lam_d}{\lam_{d,j}} \right)
< \infty.  
\end{equation}
\end{cor}

{\bf Proof.}\  Quasi-polynomial tractability implies that 
\eqref{cr_qpt} holds for some $\delta\in(0,1)$. Let
$\gamma=\gamma(d):= \frac{\delta}{\ln_+d}$. Using~\eqref{jen} we obtain
\begin{eqnarray*}
\frac{\sum_{j=1}^\infty \lam_{d,j}^{1-\frac{\delta}{\ln_+d}} }
       {\left(\sum_{j=1}^\infty \lam_{d,j}\right)^{1-\frac{\delta}{\ln_+d}}}
&=&
\frac{\sum_{j=1}^\infty \lam_{d,j}^{1-\gamma} }{\Lam_d^{1-\gamma}} 
\\
&\ge& 
\Lam_d^\gamma \ \exp\left(- \gamma  
\sum_{j=1}^\infty  \frac{\lam_{d,j} \, \ln\,\lam_{d,j}\   }{\Lam_d}\right)
\\
&=&
 \exp\left\{ \gamma \left( \ln\,\Lam_d-  \sum_{j=1}^\infty  
\frac{\lam_{d,j} \, \ln\,\lam_{d,j}}{\Lam_d} \right)  \right\}
\\
&=&
 \exp\left\{ \gamma  \sum_{j=1}^\infty
  \frac{\lam_{d,j}}{\Lam_d}\, \ln\left( \frac{\Lam_d}{\lam_{d,j}}\ 
\right) \right\}.
\end{eqnarray*}
The claim~\eqref{nes_qpt} now follows from~\eqref{cr_qpt}.
\hfill $\Box$
\bigskip

We will use later the following simple inequality that
provides a sufficient condition for the curse of dimensionality.
Recall that ${\rm trace}(C_{\nu_d})=\sum_{j=1}^\infty \lam_{d,j}$
denotes the trace of the covariance operator.

\begin{lemma} \label{cond_cd}
For any $d\in\N$ and $\eps>0$ we have 
$$
   n^{\rm avg}(\eps,d)\ge (1-\eps^2)\,\frac{{\rm trace}(C_{\nu_d})}
{\lam_{d,1}}=(1-\eps^2)\left(1+\sum_{j=2}^\infty\frac{\lam_{d,j}}{\lam_{d,1}}\right). 
$$
In particular, if ${\rm trace}(C_{\nu_d})/\lam_{d,1}\ge (1+h)^d$
for some $h>0$ and all $d\in\N$, then we have curse of dimensionality. 
\end{lemma}

{\bf Proof.}\  
For $n=n^{\rm avg}(\eps,d)$ we have
$$
   {\rm trace}(C_{\nu_d})-n\lam_{d,1} 
   \le {\rm trace}(C_{\nu_d})-\sum_{j=1}^n\lambda_{d,j}=
   \sum_{j=n+1}^\infty\lambda_{d,j}\le\eps^2\,{\rm trace}(C_{\nu_d}).
$$
Hence
$$
   n^{\rm avg}(\eps,d)\ge (1-\eps^2)\,{\rm trace}(C_{\nu_d})/\lam_{d,1},
$$
as claimed. \hfill $\Box$

\section{Tensor Products Problems} \label{s4}

In this section we assume that $F_d,G_d$ and $S_d$ are given by tensor
products. That is,
$$
F_d=F^{(1)}_{1}\otimes F^{(1)}_{2}\otimes\cdots \otimes F^{(1)}_{d}
\ \ \ \mbox{and}\ \ \
G_d=G^{(1)}_{1}\otimes G^{(1)}_{2}\otimes\cdots \otimes G^{(1)}_{d}
$$
for some Banach spaces $F^{(1)}_{k}$ of univariate real functions 
equipped with a zero-mean Gaussian measures $\mu^{(1)}_{k}$,
and some Hilbert spaces $G^{(1)}_{k}$. Here the upper index $1$ 
   reminds us that the objects are univariate.
Furthermore we assume that
$$
S_d=S^{(1)}_1\otimes S^{(1)}_2\otimes\cdots\otimes S^{(1)}_d
$$
for continuous linear operators $S^{(1)}_k:F^{(1)}_{k}\to G^{(1)}_{k}$ and
$k=1,2,\dots,d$.
 
Let $\nu^{(1)}_{k}=\mu^{(1)}_{k} (S^{(1)}_k)^{-1}$ and let 
$C^{(1)}_k :H^{(1)}_{k}\to H^{(1)}_{k}$ be 
the covariance operator of the measure $\nu^{(1)}_{k}$.
The eigenpairs of $C^{(1)}_{k}$ are denoted by
  $(\lambda(k,j),\eta(k,j))$ and 
$$
\lambda(k,1)\ge\lambda(k,2)\ge\cdots\ge 0
$$ 
as well as $\sum_{j=1}^\infty\lambda(k,j)<\infty$. 
To avoid the trivial case we assume that $\lam(k,1)>0$ for all $k\in
\N$.

The covariance operator $C_{\nu_d}$ is now the tensor product
$$
C_{\nu_d}=C^{(1)}_{1}\otimes C^{(1)}_{2}\otimes\cdots\otimes
  C^{(1)}_{d}
$$
and therefore the eigenvalues $\lambda_{d,j}$ and the eigenfunctions
$\eta_{d,j}$ are given by corresponding products of 
the one-dimensional eigenvalues and eigenvectors
$\lambda(k,j)$ and $\eta_{k,j}$, respectively. More precisely we have
$$
\{\lambda_{d,j}\}_{j\in \N}=\left\{\prod_{k=1}^d\lambda(k,j_k)
\right\}_{j_1,j_2,\dots j_d\in\N}.
$$ 
Note that
\begin{equation}\label{spps}
\sum_{j\in \N} \lam_{d,j}^\tau=\prod_{k=1}^d \sum_{j=1}^\infty \lam(k,j)^\tau
\qquad \textrm{for any}\ \tau>0.
\end{equation}
We want to express necessary and sufficient conditions, 
for each of the four types of tractability,
in terms of the eigenvalues $\lambda(k,j)$, $k,j\in\N$.
The homogeneous case of the tensor product problem, 
i.e., when $F^{(1)}_{k}=F^{(1)}_1, G^{(1)}_{k}=G^{(1)}_1$ and $S^{(1)}_k=S^{(1)}_1$ 
which implies that
$$
\lam(k,j)=\lam(1,j) \ \ \ \ \ \mbox{for all}\ \ \ k,j=1,2,\dots\,,
$$
was studied in \cite[Section 6.2] {NW08}
and in a recent paper \cite{PP}. 
In this section we mainly focus on a non-homogeneous 
case.

\subsection{Polynomial Tractability}

We know that $S=\{S_d\}$ is polynomially and strongly polynomially 
tractable iff~\eqref{poltract} holds. We now simplify 
the condition~\eqref{poltract} for tensor product problems.

\begin{thm} \label{t:tp_62} 
Consider a tensor product problem $S=\{S_d\}$.
Then
\begin{itemize}
\item 
$S$ is strongly polynomially tractable iff there exists $\tau\in (0,1)$
such that
\begin{equation}\label{cr_sptt1}
   \sum_{k=1}^\infty \sum_{j=2}^\infty 
\left(\frac{\lam(k,j)}{\lambda(k,1)}\right)^\tau  \  <\infty.
\end{equation}
If so the exponent of strong polynomial tractability is
$$
p^{\,\rm str-avg}=\inf\left\{\frac{2\tau}{1-\tau}\ \bigg|\ \ \tau\ \
  \mbox{satisfies~\eqref{cr_sptt1}}\right\}.
$$
\item
$S$ is polynomially tractable iff there exists $\tau\in (0,1)$
such that
\begin{equation}\label{cr_tp_ptln}
Q_\tau :=   \sup_{d\in \N}  \frac{1}{\ln_+d}\ \sum_{k=1}^d 
   \ln\, \left( 1+ \sum_{j=2}^\infty \left(\frac{\lam(k,j)}{\lam(k,1)}
\right)^\tau \right)  <\infty.
\end{equation}
A simpler and stronger condition
\begin{equation}\label{cr_tp_pt}
   \sup_{d\in \N} \ \frac{1}{\ln_+d}\ \sum_{k=1}^d 
   \sum_{j=2}^\infty \left(\frac{\lam(k,j)}{\lam(k,1)}\right)^\tau  \  <\infty, 
\end{equation}
is sufficient for polynomial tractability and necessary whenever
\begin{equation}
   \label{cr_tp_ptsup}
   \sup_{k\in\N}\ \sum_{j=1}^\infty \left(\frac{\lam(k,j)}{\lam(k,1)}
\right)^\tau  \  <\infty. 
\end{equation}
\end{itemize}
\end{thm}
\bigskip
{\bf Proof.}\ We prove the four conditions 
in the iff statements. Let
$$
\bb(k,j):=\frac{\lam(k,j)}{\lam(k,1)}\ \ \ \ \ 
\mbox{for all}\ \ \ k,j\in\N,
$$
be the sequence of the normalized eigenvalues so that $1=\bb(k,1)\ge
\bb(k,j)$. 
We need to verify~\eqref{poltract} which by \eqref{spps} 
now asserts that for some $q\ge0$ and $\tau\in(0,1)$ we have 
$$
C_{q,\tau}:=\sup_{d\in\N}\ \frac{\left(\sum_{j=1}^\infty\lambda_{d,j}^\tau
\right)^{1/\tau}}{\sum_{j=1}^\infty\lambda_{d,j}}\,d^{-q}=
\sup_{d\in\N}\prod_{k=1}^d\frac{\left(\sum_{j=1}^\infty
\bb(k,j)^\tau\right)^{1/\tau}}
{\sum_{j=1}^\infty\bb(k,j)}\,d^{-q}
<\infty.
$$ 
For strong polynomial tractability $q=0$, whereas for polynomial
tractability $q\ge0$.
\medskip

1. \emph{ Sufficiency of~\eqref{cr_sptt1} 
for strong polynomial tractability.} Note that 
\begin{eqnarray}
\nonumber
\prod_{k=1}^d \sum_{j=1}^\infty \bb(k,j)^\tau&=&
\prod_{k=1}^d \left( 1+ \sum_{j=2}^\infty \bb(k,j)^\tau \right) 
\le  \prod_{k=1}^d \exp\left( \sum_{j=2}^\infty \bb(k,j)^\tau \right) 
\\
\label{expbound}
&=&  \exp\left( \sum_{k=1}^d  \sum_{j=2}^\infty \bb(k,j)^\tau \right)\le
\exp\left( \sum_{k=1}^\infty \sum_{j=2}^\infty \bb(k,j)^\tau \right)<\infty
\end{eqnarray}
due to~\eqref{cr_sptt1}.
On the other hand, $\prod_{k=1}^d\sum_{j=1}^\infty\bb(k,j)\ge1$,
hence $C_{0,\tau}<\infty$. This implies strong polynomial tractability.

\medskip

2. \emph{Necessity of~\eqref{cr_sptt1} for strong polynomial tractability.}
We now know that $C_{0,\tau}<\infty$ for some $\tau\in(0,1)$.
This implies that 
$$
\prod_{k=1}^\infty\left(1+\sum_{j=2}^\infty\bb(k,j)^\tau\right)\le
C_{0,\tau}^{\,\tau} \prod_{k=1}^{\infty} \left(1+\sum_{j=2}^\infty\bb(k,j)\right)^\tau.
$$
Since $\bb(k,j)\le1$ and $\tau\in(0,1)$,
we can estimate $\bb(k,j)$ by $\bb(k,j)^\tau$. This yields
$$
\prod_{k=1}^\infty\left(1+\sum_{j=2}^\infty\bb(k,j)^\tau\right)<
C_{0.\tau}^{\,\tau/(1-\tau)}<\infty.
$$
This is equivalent to 
$$
\sum_{k=1}^\infty\sum_{j=2}^\infty\bb(k,j)^\tau=
\sum_{k=1}^\infty\sum_{j=2}^\infty
\left(\frac{\lam(k,j)}{\lam(k,1)}\right)^\tau<\infty.
$$
Hence~\eqref{cr_sptt1} holds, as claimed.
The formula for the exponent of strong polynomial tractability follows
from~\eqref{expostpol}. 
\medskip

3. \emph{Sufficiency of \eqref{cr_tp_ptln} for polynomial tractability.}
By (\ref{cr_tp_ptln}) we have
\begin{eqnarray*}
\prod_{k=1}^d \sum_{j=1}^\infty \bb(k,j)^\tau&=&
\prod_{k=1}^d \left( 1+ \sum_{j=2}^\infty \bb(k,j)^\tau \right)\\ 
&=&\exp\left(  \sum_{k=1}^d \ln \left( 1+ \sum_{j=2}^\infty 
\bb(k,j)^\tau \right) \right) 
\le \max(e^{\,Q_\tau},d^{\,Q_\tau}).
\end{eqnarray*}
Using again the fact that
$\prod_{k=1}^d\sum_{j=1}^\infty\bb(k,j)\ge1$, we conclude that
$C_{q,\tau}<\infty$ for $q=Q_{\tau}/\tau$, and  obtain polynomial tractability. 
Since condition~\eqref{cr_tp_pt}  is stronger than~\eqref{cr_tp_ptln}, 
it is also sufficient for  polynomial tractability.

\medskip

4. {\it Necessity of \eqref{cr_tp_ptln} for polynomial tractability.}
We now know that $C_{q,\tau}<\infty$ for some $q\ge 0$ and $\tau\in(0,1)$.
Proceeding as before we conclude that 
$$
\prod_{k=1}^d\left(1+\sum_{j=2}^\infty\bb(k,j)^\tau\right)
\le C_{q,\tau}^{\,\tau/(1-\tau)}\,d^{\,q\,\tau/(1-\tau)}.
$$
Hence, 
\[  
   \sum_{k=1}^d \ln \left( 1+ \sum_{j=2}^\infty \bb(k,j)^\tau\right)
  \le \frac{q\,\tau}{1-\tau}\ \ln_+d + {\frac{\tau}{1-\tau}} \ln\, C_{q,\tau},
\]
and~\eqref{cr_tp_ptln} follows.

It is easy to see that under the assumption~\eqref{cr_tp_ptsup},
the conditions~\eqref{cr_tp_ptln} and~\eqref{cr_tp_pt} are
equivalent.  Therefore,~\eqref{cr_tp_pt} is also necessary in this case.
\hfill $\Box$

\vskip 1pc

We comment on the necessary condition for polynomial
tractability. Typically, the coordinates in tensor product 
problems are ordered according to "decreasing importance". 
This means that the sequence 
$\sum_{j=2}^\infty \bb(k,j)^\tau$ is non-increasing in $k$.
In this case~\eqref{cr_tp_ptsup} holds and 
the simple condition~\eqref{cr_tp_pt} is necessary and sufficient
for polynomial tractability. 
However, in general, nothing prevents us from a strange ordering 
of important and  unimportant coordinates so that the sequence 
of $\sum_{j=2}^\infty \bb(k,j)^\tau$ is not non-increasing in $k$.
In this case the stronger condition~\eqref{cr_tp_pt} may fail as 
illustrated by the following example. Let      
\[
  \lam(k,j)=\bb(k,j)=\begin{cases}
  1 &\ \ \ \mbox{for}\ \ \  j=1,\\
  1 &\ \ \ \mbox{for}\ \ \ j\in[2,k] \ \mbox{and}\ k=2^{2^m}
 \  
\mbox{for non-negative  integer $m$},\\
  0 & \ \ \ \textrm{otherwise}.
  \end{cases}
\]
By counting the number of $1$'s in $\lam(k,j)$ we easily conclude that 
$$
n^{\rm avg}(\eps,d)\le \prod_{m\in\N,\,2^{2^m}\le d}2^{2^m}
\le 2^{2^{\ln_2(\ln_2(d^{\,2}))}}=d^{\,2}\ \ \ \ \
\mbox{for all}\ \ \ \eps\in[0,1)\ \ \mbox{and}\ \ d\in \N.
$$
So polynomial tractability holds but
condition~\eqref{cr_tp_pt} fails. Therefore, in general,
it is not necessary for polynomial tractability.

\subsection{Quasi-Polynomial Tractability}

We now consider quasi-polynomial tractability of tensor products.
First of all let us check how the right-hand side of
Lemma~\ref{useful1} simplifies in this
case. 
Let 
$$
\Lambda(k):=\sum_{j=1}^\infty\lam(k,j)\ \ \ \mbox{and by}\  \eqref{spps}\ \ \
\Lambda_d:=\sum_{j=1}^\infty\lambda_{d,j}=\prod_{k=1}^d\Lambda(k).
$$
For tensor products we have 
\begin{eqnarray} \nonumber
 \sum_{j=1}^\infty  \lam_{d,j} \, \ln\,\lam_{d,j}\
&=& \nonumber
 \sum_{z=[z_1,z_2,\dots,z_d]\in \N^d}
\  \prod_{k=1}^d \lam(k,z_k) \, 
\sum_{k=1}^d \ln\,\lam(k,z_k)
\\
&=& \nonumber
\sum_{k=1}^d \sum_{z\in \N^d}  \lam(k,z_k) \, \ln\,\lam(k,z_k)   
\prod_{{1\le m\le d\atop m\not = k}} \lam(m,z_m) 
\\ \nonumber
&=& 
\sum_{k=1}^d \left( \sum_{j=1}^\infty  \lam(k,j) \, \ln\,\lam(k,j) \right)   
\prod_{{1\le m\le d\atop m\not = k}} \left( \sum_{j=1}^\infty \lam(m,j) \right)
\\ \nonumber
&=& 
\sum_{k=1}^d \left( \sum_{j=1}^\infty  \lam(k,j) \, \ln\,\lam(k,j) \right)   
\prod_{{1\le m\le d\atop m\not = k}} \Lam(m)
 \\
&=&  \label{tp_ln}
\sum_{k=1}^d \left( \sum_{j=1}^\infty  \lam(k,j) \, \ln\,\lam(k,j) \right)   
\frac{\Lam_d}{\Lam(k)}.
\end{eqnarray}
Inequality~\eqref{jen} now becomes
\begin{equation}
\label{jent}
 \Lam_d^{-1} \sum_{j=1}^\infty\lam_{d,j}^{1-\g} 
\ge
\exp\left( - \g  
\sum_{k=1}^d \frac1{\Lambda(k)}\,\sum_{j=1}^\infty 
\lam(k,j) \, \ln\,\lam(k,j) \right).
\end{equation}
This inequality will be used in the following theorem
which addresses quasi-polynomial tractability for tensor product
problems.

\begin{thm} \label{t:tp_qpt} 
Consider a tensor product problem $S=\{S_d\}$. Then
\begin{itemize}
\item $S$ is quasi-polynomially tractable iff there exists $\delta\in(0,1)$
such that
\begin{equation}\label{cr_tp_qpt}
   \sup_{d\in\N}\  
\prod_{k=1}^d \frac{\sum_{j=1}^\infty \lam(k,j)^{1-\frac{\delta}{\ln_+d}}}
         {\left(\sum_{j=1}^\infty\lam(k,j)\right)^{1-\frac{\delta}{\ln_+d}}}
  <\infty.
\end{equation}
\item If $S$ is quasi-polynomially tractable then 
\begin{equation}\label{nes_tp_qpt}
\sup_{d\in \N} \ \frac{1}{\ln_+d} \ 
\sum_{k=1}^d  \sum_{j=1}^\infty  \frac{\lam(k,j)}{\Lam(k)} \, 
\ln\left(\frac{\Lam(k)} {\lam(k,j)} \right) 
<\infty.
\end{equation}
\item If there exists $\delta>0$ such that
\begin{equation}\label{cr_tp_qptln}
   \sup_{d\in \N}\  \sum_{k=1}^d 
   \ln\left(1+\sum_{j=2}^\infty \left(\frac{\lam(k,j)}
{\lam(k,1)}\right)^{1-\frac{\delta}{\ln d}}  
\right) <\infty, 
\end{equation}
or
\begin{equation}\label{cr_tp_qptsum}
   \sup_{d\in \N} \ \sum_{k=1}^d 
   \sum_{j=2}^\infty \left(\frac{\lam(k,j)}{\lam(k,1)}
\right)^{1-\frac{\delta}{\ln d}}  \  <\infty
\end{equation}
then $S$ is quasi-polynomially tractable.
\end{itemize}
\end{thm}

{\bf Proof.} 
In view of  \eqref{spps}, 
criterion~\eqref{cr_tp_qpt} is just the general criterion~\eqref{cr_qpt} 
in Theorem \ref{t:qpt}  specified for tensor products.
The necessary condition in~\eqref{nes_tp_qpt} 
is just a specification of the general necessary
condition in~\eqref{nes_qpt} for tensor products. To see this,
note that 
\begin{eqnarray*}
\sum_{j=1}^\infty  \frac{\lam_{d,j}}{\Lam_d} \, 
\ln \left(\frac {\Lam_d} {\lam_{d,j}}\right)
&=& 
\ln\,\Lam_d - \Lam_d^{-1} \sum_{j=1}^\infty  \lam_{d,j} \, \ln\, \lam_{d,j}
\\
(\textrm{by}\ (\ref{tp_ln}) )  
&=&
 \sum_{k=1}^d \ln\, \Lam(k) - 
\sum_{k=1}^d \sum_{j=1}^\infty \frac{\lam(k,j)}{\Lam(k)} \, \ln\,\lam(k,j) 
\\  
&=& 
\sum_{k=1}^d \sum_{j=1}^\infty  \frac{\lam(k,j)}{\Lam(k)} \, 
\ln\left(\frac{\Lam(k)} {\lam(k,j)} \right) .
\end{eqnarray*}

To see that~\eqref{cr_tp_qptln} 
is sufficient for quasi-polynomial tractability,
observe that the fraction in~\eqref{cr_tp_qpt} can be written
with $\bb(k,j)=\lam(k,j)/\lam(k,1)$  as
$$
\prod_{k=1}^d\,\frac{1+\sum_{j=2}^\infty\bb(k,j)^{1-\frac{\delta}{\ln_+d}}}
{\left(1+\sum_{j=2}^\infty\bb(k,j)\right)^{1-\frac{\delta}{\ln_+d}}}.
$$
Taking logarithms, we see that the numerator is bounded
by~\eqref{cr_tp_qptln} while the denominator is larger than 1. 
Hence~\eqref{cr_tp_qpt} is bounded and we are done.

Since (\ref{cr_tp_qptsum}) is stronger 
than (\ref{cr_tp_qptln}), it is also sufficient for 
quasi-polynomial tractability. $\Box$
\bigskip

\subsection{Weak Tractability}

We present a simple criterion of weak tractability for tensor products.

\begin{thm} \label{t:tp_wt} 
Consider a tensor product problem $S=\{S_d\}$. 
If for some $\tau\in (0,1)$
\begin{equation}\label{cr_tp_wt}
   \lim_{d\to\infty}\,  d^{-1}\, \sum_{k=1}^d \sum_{j=2}^\infty 
\left(\frac{\lam(k,j)}{\lam(k,1)}\right)^{\tau}=0
\end{equation}
then $S$ is weakly tractable.
\end{thm}
{\bf Proof.}
The idea is basically the same as in the proof of Theorem \ref{t:tp_62}.
Namely, we apply~\eqref{cheb} with $z=1$. 
As before, let $\bb(k,j):=\lam(k,j)/\lam(k,1)$.  Then~\eqref{cheb},
by \eqref{spps} can be rewritten as
$$
n^{\rm avg}(\eps,d)\le
\prod_{k=1}^d\left[\frac{1+\sum_{j=2}^\infty\bb(k,j)^\tau}{
1+\sum_{j=2}^\infty\bb(k,j)}\right]^{1/(1-\tau)}\,\eps^{-2/(1-\tau)}.
$$
Since the denominator above is larger than one, it may be dropped. 
Using~\eqref{expbound} we have 
\[
   n^{\rm avg}(\eps,d)\le \left[ 
    \exp\left( \sum_{k=1}^d  \sum_{j=2}^\infty \bb(k,j)^\tau 
\right) \eps^{-2} \right]^{1/(1-\tau)}
    = \exp\left[ (1-\tau)^{-1} (\theta_d d + 2 \ln\, \eps^{-1})\right],   
\]
where
\begin{equation} \label{thd}
  \theta_d:= d^{-1} \sum_{k=1}^d  \sum_{j=2}^\infty \bb(k,j)^\tau 
  \to 0 \qquad \textrm{as}\ \ d\to \infty
\end{equation}
due to~\eqref{cr_tp_wt}. Equivalently,
 \[
    \ln\,n^{\rm avg}a(\eps,d) \le (1-\tau)^{-1} 
\left[ \theta_d\, d + 2 \ln\, \eps^{-1}\right].
 \] 
By~\eqref{thd}
\[
  \lim_{d+\eps^{-1}\to \infty}  
  \frac
  {\theta_d \, d + 2 \ln\, \eps^{-1}}
  {d+ \eps^{-1}}=0,
\]
and we obtain the weak tractability. \hfill $\Box$

\vskip 1pc
Note that~\eqref{cr_tp_wt} holds if
\begin{equation}\label{cr_tp_wt2}
   \lim_{k\to\infty}   \sum_{j=2}^\infty \left(\frac{\lam(k,j)}
{\lam(k,1)}\right)^{\tau}=0.
\end{equation}
Hence~\eqref{cr_tp_wt2} implies weak tractability. The last condition
yields
$$
\lim_{k\to\infty}\frac{{\rm trace}(C^{(1)}_{k})}{\lam(k,1)}=1
$$ 
so that the Gaussian measure is asymptotically 
concentrated on the one-dimensional subspace ${\rm span}(\eta(k,1))$ of
$H^{(1)}_{k}$.

\section{Multivariate Approximation and Korobov Kernels} \label{s5}

The non-homogeneous case offers the possibility of vanquishing the
curse of dimensionality via variation of weights and smoothness parameters. 
We illustrate this by an example with  Korobov kernels of 
decreasing weights $g_k$ and increasing smoothness $r_k$. 
As we shall see, even strong polynomial 
tractability holds if the decay of $g_k$ is sufficiently fast.
Multivariate approximation for Korobov spaces 
in the worst case setting was recently
studied  in~\cite{PW10}. Here we present its average case analog.

In this section we consider a multivariate approximation problem
defined over the space of continuous real functions
equipped with a zero-mean Gaussian measure whose covariance is given as 
a Korobov kernel. More precisely, consider the approximation problem 
$${\rm APP}=\{ {\rm APP}_d \}_{d\in\N}\ \ \  \mbox{with}\ \ \ 
{\rm APP}_d: C([0,1]^d)\to L_2([0,1]^d) 
$$
given by
$$
{\rm APP}_df=f \ \ \ \mbox{for all} \ \ \ f\in C([0,1]^d).
$$
The space  $C([0,1]^d)$ of continuous real functions
is equipped with a zero-mean Gaussian measure~$\mu_d$ whose 
covariance kernel
$$
K_d(x,y)= \int_{C([0,1]^d)} f(x)f(y)\, \mu_d(df), \quad
x,y\in[0,1]^d,
$$
is given as follows. First of all we assume that $K_d$ is of product
form, 
$$
K_d(x,y)= \prod_{k=1}^d \RR_k(x_k,y_k)\ \ \ \ \ 
\mbox{for all}\ \ \ x=[x_1,x_2,\dots,x_d],\ y=[y_1,y_2,\dots,y_d]\in[0,1]^d,
$$
where $\RR_k=\RR_{r_k,g_k}$ are univariate Korobov kernels,
$$
   \RR_{r,\beta}(x,y):=1+2\,\beta\,
   \sum_{j=1}^\infty  j^{-2r}\,\cos(2\pi\,j(x-y))\ \ \ \ \
\mbox{for all}\ \ \ x,y\in[0,1].
$$ 
Here $\beta\in (0,1]$ and $r$ is a real number such that 
$r>\tfrac12$. Note that for $y=x$ we have
$$
\RR_{r,\beta}(x,x)=1+2\,\beta\, \,\zeta(2r),
$$
where $\zeta(x)=\sum_{j=1}^\infty j^{-x}$ is the Riemann zeta
function which is well-defined only for $x>1$. 
That is why we have to consider $r>\tfrac12$. 

We assume that the sequence $\{r_k\}$ is non-decreasing,
\begin{equation} \label{rk}
  \tfrac12 < r_1\le r_2\le\cdots\le r_d\le\cdots.
\end{equation}

The weight sequence $\{g_k\}$ serves as a scaling and, as we shall see,
tractability results will depend on the behavior of $g_k$ at infinity. 
We assume that
\begin{equation} \label{gk}
  1 \ge g_1\ge g_2 \ge\cdots >0.
\end{equation}

As already mentioned, the sequences $\{r_k\}$ and $\{g_k\}$ may be related,
$g_k=g(r_k)$ for some non-increasing function $g:[\tfrac12,\infty)\to[0,1]$.
The case which can be often found in the literature corresponds
to $g_k=1$ or $g_k=(2\pi)^{-2r_k}$. For $g_k=g(r_k)$ the behavior of $g_k$ 
at infinity depends on the function $g$ and the behavior of $r_k$ at infinity. 
A summary of the properties of the Korobov kernels can be found 
in Appendix A of~\cite{NW08}.

For a fixed $d$, the multivariate approximation problem
under similar conditions was studied in \cite{MW,PGW}. 
For varying $d$, the homogeneous case, i.e., $\RR_k=\RR$
for all $k$  with $\RR$ not necessarily equal to a Korobov kernel, 
was studied in~\cite{LT,LZ,NW08}. In this case, 
we have the curse of dimensionality since $n^{\rm avg}(\eps,d)$
depends exponentially on $d$. 

The induced measure $\nu_d=\mu_d {\rm APP}_d^{-1}$ on 
$L_2([0,1]^d)$ is also a zero-mean Gaussian
measure. It is known, see e.g., \cite{NW08}, that 
the eigenvalues of its covariance operator $C_{\nu_d}$ are
given by 
\begin{equation}\label{eq:korevals1}
   \lambda_{d,z}=\prod_{k=1}^d \lambda(k,z_k) \ \ \ \ \ \mbox{for all}\ \ \ 
   z=[z_1,z_2,\dots,z_d]\in\N^d,
\end{equation}
where $\lambda(k,1)=1$ and 
\begin{equation}\label{eq:korevals2}
    \lambda(k, 2j)= \lambda(k, 2j+1)= \frac{g_k}{j^{2r_k}}\ , \qquad j\in \N.
\end{equation}
Note that the trace of $C_{\nu_d}$ is
$$
{\rm trace}(C_{\nu_d})=\prod_{k=1}^d
\left(1+2\,g_k\,\zeta(2r_k)\right).
$$

We have the curse of dimensionality when 
$$
   g_{\rm \,lim}:=\lim_{k\to\infty}\,g_k>0.
$$
Indeed, in this case 
$$
   {\rm trace}(C_{\nu_d})\ge (1+2\,g_{\rm \,lim})^d,
$$
and Lemma \ref{cond_cd} yields the curse.
Therefore $\lim_k g_k=0$ is a necessary condition to vanquish 
the curse.

\begin{thm} \label{t:th_korob}
Consider the approximation problem ${\rm APP}=\{ {\rm APP}_d\}$ 
in the average case with a zero-mean Gaussian measure 
whose covariance operator is given as the Korobov kernel with
the weights $g_k$ and smoothness $r_k$ satisfying 
\eqref{gk} and \eqref{rk}, respectively. Then 
\begin{itemize}
\item ${\rm APP}$ is polynomially tractable iff
\begin{equation}\label{pt_kor}
   \rho_g:=\liminf_{k\to\infty}\ \frac{\ln\, \frac1{g_k}}{\ln k} >1.
\end{equation}
\item ${\rm APP}$ is strongly polynomially tractable 
iff it is polynomially tractable. 
If so, the exponent of strong polynomial tractability is 
$$
   p^{\rm avg-str}= \max\left( \frac{2}{2r_1-1}, 
   \frac{2}{\rho_g -1}\right).
$$
\item
If ${\rm APP}$ is quasi-polynomially tractable then  
\begin{equation}\label{qpt_kor}
    \sup_{d\in\N}\ \frac1{\ln_+d}\ \sum_{k=1}^d g_k\,
    \ln_+\,\frac1{g_k}<\infty.
\end{equation}
If~\eqref{qpt_kor} holds and
\begin{equation}\label{needed?}
     \liminf_{k\to\infty}\frac{r_k}{\ln k}>0
\end{equation}
then ${\rm APP}$ is quasi-polynomially tractable.

\item ${\rm APP}$ is weakly tractable iff
$$
    \lim_{k\to \infty} g_k=0.
$$
\end{itemize}

\end{thm}

\noindent{\bf Proof:} 
We will use Theorem \ref{t:tp_62} and proceed in a way similar to that
of the proof of Theorem~1 in~\cite{PW10}. 
The main difference is that here $\tau\in(0,1)$.

We first show that \eqref{pt_kor} implies 
strong polynomial tractability and then that
polynomial tractability implies \eqref{pt_kor}. 
Assume thus that \eqref{pt_kor} is satisfied.
Then for some $\delta>0$ and  all large $k$ we have 
\[
   \frac{\ln\,\frac1{g_k}}{\ln\,k}\ge 1+\delta.
\] 
Hence, there is a positive $C$ such that for any $\tau\in (0,1)$ we
obtain 
$$
   g_k^\tau\le C^{\tau}\,k^{-\tau(1+\delta)}\ \ \ \ \ 
   \mbox{for all}\ \ \ \ \  k\in \N.
$$
If we choose $\tau\in (\tfrac{1}{1+\delta},1)\cap (\tfrac1{2r_1},1)$ then
\begin{eqnarray*}
  \sum_{k=1}^\infty \sum_{j=2}^\infty 
  \left(\frac{\lam(k,j)}{\lam(k,1)}\right)^\tau 
  &=& 2 \sum_{k=1}^\infty g_k^\tau\,\sum_{j=1}^\infty j^{-2\tau r_k}
  \le   2C^{\tau} \sup_k  \zeta(2\tau r_k)  \sum_{k=1}^\infty  k^{-\tau(1+\delta)}
\\
  &\le&  2C^{\tau}  \zeta(2\tau r_1)\,\zeta\left(\tau(1+\delta)\right)<\infty,
\end{eqnarray*}
and condition \eqref{cr_sptt1} 
of Theorem \ref{t:tp_62} yields strong polynomial tractability.

Assume now that polynomial tractability holds. 
Then for $\tau\in(\tfrac1{2r_1},1)$ we have 
\[ 
  \sup_k \sum_{j=2}^\infty \left(\frac{\lam(k,j)}{\lam(k,1)}\right)^\tau 
  =2 \sup_k  [g_k^\tau\,\zeta(2\tau r_k)]
  =2 g_1^\tau\,\zeta(2\tau r_1) <\infty.
\]
Therefore, condition  \eqref{cr_tp_ptsup} is verified, 
hence condition  \eqref{cr_tp_pt}
is necessary for polynomial tractability.
The latter condition for the Korobov case is 
\[
   C:= \sup_{d\in \N} \ \frac{2}{\ln_+d}\ 
        \sum_{k=1}^d g_k^\tau\,\zeta(2\tau r_k) \  <\infty, 
\]
for some $\tau\in(\tfrac1{2r_1},1)$.
All terms of the last sum are larger or equal to $g_d^\tau$ and therefore
for $d>1$ we have  $g_d^\tau\le \tfrac{C\ln d}{2d}$.
This is equivalent to
\[
   \frac{\ln\,\frac1{g_d}}{\ln d} \ge \frac{1}{\tau}
   \left(1-\frac{\ln(C/2)+\ln\ln d}{\ln d}\right).
\]
Hence,
\[
  \rho_g= \liminf_{d\to\infty} \frac{\ln\,\frac1{g_d}}{\ln d} 
   \ge  \frac{1}{\tau}  > 1,
\]
as required in  \eqref{pt_kor}.

We now turn to the exponent of strong polynomial tractability. We must have 
$\tau>\tfrac{1}{2 r_1}$ and from the last displayed formula
$\tau> \tfrac{1}{\rho_g}$. 
>From Theorem~\ref{t:tp_62} we obtain that
$$p^{{\rm avg-str}}=
\max\left( \frac 2{2r_1 -1}, \frac{2}{\rho_g -1}\right).
$$
This completes the proof of polynomial tractability.
\vskip 1pc

Assume now that quasi-polynomial tractability  holds. 
Then the necessary condition \eqref{nes_tp_qpt} is satisfied. 
Clearly, all terms appearing in this condition are
positive. We simplify  \eqref{nes_tp_qpt} by omitting all terms for
$j\not=2$, and obtain
\begin{equation}\label{conse1}
   \sup_{d\ge\N} \ \frac{1}{\ln_+d} \  
   \sum_{k=1}^d  \frac{\lam(k,2)}{\Lam(k)} \, 
   \ln\left(\frac{\Lam(k)} {\lam(k,2)} \right) <\infty.
\end{equation}
Recall that for the Korobov case,
$\Lam(k)=1+2\,g_k\,\zeta(2r_k)$ and $\lam(k,2)=g_k$. 
Since $\Lam(k)\ge 1$ 
and $\Lam(k)/\lam(k,2)\ge3$ 
we obtain
\[
   \sup_{d\in\N} \ \frac{1}{\ln_+d} \ 
   \sum_{k=1}^d  \frac{\lam(k,2)}{\Lam(k)} 
   \, \ln_+\left(\frac{1} {\lam(k,2)} \right) <\infty.
\]
Furthermore, since $\{\Lam(k)\}$ is non-increasing, we have 
\[
   \sup_{d\in\N} \ \frac{1}{\ln_+d} \ 
   \sum_{k=1}^d  \lam(k,2) \, \ln_+\left(\frac{1} {\lam(k,2)} \right) 
   <\infty.
\]
This is equivalent to~\eqref{qpt_kor}, and completes this part of the proof. 

We now prove that~\eqref{qpt_kor} and~\eqref{needed?} are sufficient 
for quasi-polynomial tractability. 
Theorem~\ref{t:tp_qpt}  states that ${\rm APP}$ is quasi-polynomially
tractable iff there exists $\delta\in(0,1)$ such that~\eqref{cr_tp_qpt}
holds, i.e., 
\begin{equation}\label{quasicondition}
   \sup_{d\in\N}\ \prod_{k=1}^d 
   \frac{ \sum_{j=1}^\infty \lambda(k,j)^{\tau_d}}
  {\left(  \sum_{j=1}^\infty\lambda(k,j)\right)^{\tau_d}} <\infty,
\end{equation}
where $\tau_d=1-\tfrac{\delta}{\ln_+d}$. 
Take any $\delta\in(0,\min(\tfrac12,1-\tfrac1{2r_1}))$. 
Inequality $\delta<1-1/(2r_1)$ ensures that all the
sums above are finite because
$2r_k\tau_d \ge 2r_1\tau_1>1$. 
 
We split the product in \eqref{quasicondition} into two products
\[
  \Pi_1(d) := \prod_{k=1}^d 
  \left(\sum_{j=1}^\infty \lambda(k,j)\right)^{\frac{\delta}{\ln_+d}}  
\]
and
\[
   \Pi_2(d):=  \prod_{k=1}^d 
            \frac{ \sum_{j=1}^\infty \lam(k,j)^{\tau_d}}
            { \sum_{j=1}^\infty \lam(k,j)  }.
\]
In what follows we will write $C$ for some positive number which is 
independent of $d$ and~$k$, and whose value may change for  
successive estimates. 

For $\Pi_1(d)$ we use $(1+x)^t=\exp(t\ln(1+x))\le
\exp(tx)$ and have 
\begin{eqnarray*}
  \Pi_1(d) &=& \prod_{k=1}^d 
  \left(  1+ \sum_{j=2}^\infty \lam(k,j)   \right)^{\frac{\delta}{\ln_+d}}  
  \le \exp\left( \frac{\delta}{\ln_+d}   \sum_{k=1}^d  
  \sum_{j=2}^\infty \lam(k,j)  \right) 
  \\
  &\le&   \exp\left( \frac C{\ln_+d}   
\sum_{k=1}^d g_k\,\zeta(2r_k) \right)\le
\exp\left( \frac{C\,\zeta(2r_1)}{\ln_+d}   
\sum_{k=1}^d g_k \right). 
\end{eqnarray*}
Clearly, \eqref{qpt_kor} implies that $\sup_{d\in\N} \Pi_1(d) <\infty$.

We now turn to the product $\Pi_2(d)$. We estimate each of its factors by 
\begin{equation} \label{suf0} 
  \frac{ \sum_{j=1}^\infty \lam(k,j)^{\tau_d}}
            { \sum_{k=1}^\infty \lam(k,j)     } 
  \le
\frac{ 1+ 2\lam(k,2)^{\tau_d}}{ 1+  2\lam(k,2)}  
  + \sum_{j=4}^\infty \lam(k,j)^{\tau_d} .            
\end{equation}

Note that if $|\ln\lam(2,k)|\le 3\ln_+d$, then
\begin{eqnarray*}
 &&  \frac{ 1+ 2\lam(k,2)^{\tau_d}}{ 1+  2\lam(k,2)} 
   =  \frac{ 1+2\lam(k,2) 
   \exp\left({\frac{-\delta\,\ln\lam(k,2)}{\ln_+d}}\right) } { 1+2\lam(k,2)} 
 \\
   &\le& \frac{ 1+2\lam(k,2)\left({1+ \frac{C|\ln\lam(k,2)|}{\ln_+d}}\right)}
   {1+  2\lam(k,2)} 
   \le  1+  \frac{C\lam(k,2) |\ln\lam(k,2)|} {\ln_+d}\, ,
\end{eqnarray*}
while if $|\ln \lam(k,2)|\ge 3\ln_+d$ then $\delta<\tfrac12$ implies
\[
   \frac{ 1+ 2\lam(k,2)^{\tau_d} } { 1+  \lam(k,2)  }
   \le  1+ 2\lam(k,2)^{\tau_d} 
   \le 1+ 2\lam(k,2)^{1/2} 
   \le  1+ 2d^{-3/2}.
\]
Thus, in any case
\begin{equation}\label{suf1}
   \frac{ 1+ 2\lam(k,2)^{\tau_d} } { 1+  2\lam(k,2)     } 
   \le 1+ 2d^{-3/2} + \frac{C\lam(k,2) |\ln\lam(k,2)|} {\ln_+d}.
\end{equation}

It remains to evaluate the sum in \eqref{suf0}. 
An easy and elementary calculation 
shows that \eqref{qpt_kor} implies $\lam(k,2)=g_k \le \tfrac Ck$.
On the other hand, \eqref{needed?} 
yields $r_k\ge h\ln k-C$ for all $k\in N$ with appropriate
$h,C>0$. We obtain now
\begin{eqnarray}\nonumber 
   \sum_{j=4}^\infty \lam(k,j)^{\tau_d}
   &\le&  C \,\lam(k,4)^{\tau_d}
    =  C \, \lam(k,2)^{\tau_d} 4^{-r_k \tau_d}
   \\ \label{suf2}
    &\le& C\cdot (C/k)^{1-\delta/\ln_+ d} \,2^{- (h\ln k-C)} \le C k^{-(1+u)},
\end{eqnarray}
where 
$u=h\ln 2>0$.
Combining~\eqref{suf0},~\eqref{suf1}) and~\eqref{suf2},
and using again $1+x\le \exp(x)$ we easily check that 
\[
    \frac{ \sum_{j=1}^\infty \lam(k,j)^{\tau_d} }
            { \sum_{j=1}^\infty \lam(k,j)     } 
    \le \exp\left(  2d^{-3/2} + \frac{C\,\lam(k,2) |\ln\lam(k,2)|} {\ln_+d} 
        + C  k^{-(1+u)}  \right).
\]
Then it follows that
\begin{eqnarray*}
   \Pi_2(d) &\le& \exp\left ( \sum_{k=1}^d \left(  2d^{-3/2} + 
   \frac{C\,\lam(k,2) |\ln\lam(k,2)|} {\ln_+d} 
    + C k^{-(1+u) } \right) \right)
   \\
   &\le& \exp\left(\sum_{k=1}^d \left(  2d^{-3/2} 
     + \frac{C\,g_k\,\ln_+\frac1{g_k}}{\ln_+d}
     + C  k^{-(1+u)}   \right)\right),
\end{eqnarray*}
and~\eqref{qpt_kor} implies that  $\sup_{d\in\N} \Pi_2(d) <\infty$. 
Therefore, 
 \[
    \sup_{d\in\N} \Pi_1(d)\,\Pi_2(d) \le  \sup_{d\in\N} \Pi_1(d) \,
    \ \sup_{d\in\N} \Pi_2(d) <\infty.
\] 
Hence,~\eqref{quasicondition} holds so that 
the quasi-polynomial tractability is proved.
\medskip

We now consider weak tractability.
Sufficiency. Let $\lim_kg_k=0$. Then for an arbitrarily small
positive $\delta$
there exists $k(\delta)$ such that $g_k\le \delta$ for all $k\ge k(\delta)$.
We check the assumption~\eqref{cr_tp_wt} of Theorem~\ref{t:tp_wt}.    
For $\tau\in(1/(2r_1),1)$ and $d>k(\delta)$ we have 
\begin{eqnarray*}
a_d:=\frac 1d\ \sum_{k=1}^d
\sum_{j=2}^\infty\left( \frac{\lam(k,j)}{\lam(k,1)}\right)^\tau
&=&
\frac2d\,\sum_{k=1}^d g_k^\tau \, \zeta(2r_k\tau)
\\
&\le& \frac{2\zeta(2r_1\tau)\,k(\delta)}{d}+ \frac{(d-k(\delta))\,\delta}{d}.
\end{eqnarray*}
Hence,
$$
\limsup_{d\to\infty}\,a_d\le \delta.
$$
{}For $\delta$ tending to zero, we conclude that 
$\limsup_d\,a_d=\lim_d\,a_d=0$,
and obtain weak tractability due to Theorem~\ref{t:tp_wt}.    

Necessity. We have already showed that $\lim_kg_k=0$
is a necessary condition for weak tractability.
This completes the proof. \hfill $\Box$
\medskip

We do not know if~\eqref{needed?} is needed for quasi-polynomial tractability. 
However, for $g_k=g(r_k)$ 
with $g(r)=\ttt^{\,r}$ and $\ttt\in(0,1)$, or
$g(r)=r^{-s}$ and  $s>0$, 
this condition clearly follows from~\eqref{qpt_kor}
since the latter implies that $g_k\le \tfrac{C}{k}$. 
For such weights and smoothness parameters,~\eqref{qpt_kor} is
a necessary and sufficient condition for quasi-polynomial tractability.
 
We illustrate Theorem~\ref{t:th_korob} for special weights.
\begin{itemize}
\item 
Let $g_k=v^{r_k}$ with $v\in(0,1)$. 
\begin{itemize}
\item Strong polynomial tractability holds $\quad\mbox{iff}\quad
\rho_r:=\liminf_{k\to\infty}\frac{r_k}{\ln\,k}>\frac{1}{\ln\,v^{-1}}$.

If so the exponent is
$\quad
p^{{\rm avg-str}}=\max\left(\frac2{2r_1-1},\frac2{\rho_r[\ln\,v^{-1}]-1}\right).
$ 
\item Quasi-polynomial tractability holds $\quad\mbox{iff}\quad
\sum_{k=1}^dv^{r_k}\,\max(1,r_k)=\mathcal{O}(\ln d)$.
\item Weak tractability holds $\quad\mbox{iff}\quad \lim_{k\to\infty}r_k=\infty.$
\end{itemize}
\item Let $g_k=r_k^{-s}$ for $s>0$. 
\begin{itemize}
\item Strong polynomial tractability holds $\quad\mbox{iff}\quad
\rho_r:=\liminf_{k\to\infty}\frac{\ln r_k}{\ln\,k}>\frac{1}{s}$.

If so the exponent is
$\quad
p^{{\rm avg-str}}=\max\left(\frac2{2r_1-1},\frac2{\rho_r\,s-1}\right).
$
\item Quasi-polynomial tractability holds $\quad\mbox{iff}\quad
\sum_{k=1}^dr_k^{-s}\,\max(1,\ln\,r_k)=\mathcal{O}(\ln d)$.
\item Weak tractability holds $\quad\mbox{iff}\quad
\lim_{k\to\infty}r_k=\infty$. 
\end{itemize}
\end{itemize}

It is also important to notice that Theorem~\ref{t:th_korob} 
holds for constant smoothness parameters $r_k\equiv r>\tfrac12$
if $g_k$ are \emph{not} related to $r_k$ and satisfy the conditions
presented in Theorem~\ref{t:th_korob}. 
This corresponds to appropriately decaying product weights, the case
that was also studied in~\cite{NW08} p. 276.

\section{Comparison of Korobov, Euler, and Wiener Kernels} \label{s6}

Another application of our general results is given in \cite{LPW1},
where tensor products of multi-parametric 
Wiener and Euler integrated processes are considered.
We briefly summarize the results of \cite{LPW1} to compare them 
to the results of the previous section.

Let $W(t), t\in[0,1]$, be a standard Wiener process, i.e. a Gaussian random
process with zero mean and covariance
$
  K_{1,0}^\ee(s,t)=K_{1,0}^\ww(s,t):=\min(s,t).
$
Consider two sequences of integrated random processes 
$X^\ww_r, X^\ee_r$ on $[0,1]$ 
defined inductively on $r$ by $X^\ww_0=X^\ee_0=W$, and for $r=0,1,2,\dots$ 
\begin{eqnarray*}
  X^\ww_{r+1}(t)&=& \int_0^t X^\ww_{r}(s) {\rm d} s,  
  \\
  X^\ee_{r+1}(t)&=& \int_{1-t}^1 X^\ee_{r}(s) {\rm d} s. 
\end{eqnarray*}
$\{X^\ww_{r}\}$ is called the 
univariate integrated Wiener process, while 
$\{X^\ee_r\}$ is called the univariate integrated Euler
process.

Clearly, $X^\ww_r$ and $X^\ee_r$ have the same smoothness
properties but they satisfy different boundary conditions.

The covariance kernel of $X^\ww_{r}$ is given by
$$
  K_{1,r}^{\ww}(x,y)=\int_0^{\min(x,y)} \frac{(x-u)^r}{r!}\,
  \frac{(y-u)^r}{r!}\,{\rm d}u  
$$
and is called the Wiener kernel, 
while the covariance kernel of $X^\ee_{r}$ is given by
$$
   K_{1,r}^{\ee}(x,y)=
   \int_{[0,1]^r}\min(x,s_1)\,\min(s_1,s_2)\,\cdots\,\min(s_r,y)\,{\rm d}s_1
   \,{\rm d}s_2\cdots{\rm d}s_r
$$
and is called the Euler kernel.
The last kernel can be expressed in terms of Euler polynomials, 
hence the name of the process and its kernel.

The corresponding tensor product kernels on $[0,1]^d$ are given by
$$
  K_d^{\ww}(s,t)=\prod_{k=1}^d K_{1,r_k}^{\ee}(s_k,t_k),\ \ \ \ \
\mbox{and}\ \ \ \ \  
   K_d^{\ee}(s,t)=\prod_{k=1}^d K_{1,r_k}^{\ee}(s_k,t_k).
$$
As for the Korobov case, the sequence $\{r_k\}$ with integers
$$
r_1\le r_2\le\cdots \le r_d\le\cdots,
$$ 
describes the increasing smoothness of a process
with respect to 
the successive coordinates.  

We now compare tractability results for processes described by the 
Euler, Korobov and Wiener kernels
from~\cite{LPW1} and from Theorem \ref{t:th_korob}. 
Some results are the same:
\begin{itemize}
\item 
strong polynomial tractability and polynomial tractability are equivalent,
\item there is a lim-inf-type criterion for polynomial tractability,
\item there is a narrow zone where quasi-polynomial tractability holds while 
 polynomial tractability fails,
\item weak tractability is equivalent to a convergence without rate, 
$\lim_k r_k=\infty$ for both integrated processes, or to 
$\lim_k g_k=0$ for Korobov case,  
\item if weak tractability fails then the curse of dimensionality appears.
\end{itemize}

The conditions on strong polynomial tractability 
for Euler and Wiener integrated processes
are different. Namely, strong polynomial tractability holds iff
\begin{eqnarray*}
   \rho_E&:=&\liminf_{d\to\infty}\ \frac{r_d}{\ln d} \ >\ \frac{1}{2\ln\, 3} 
   \qquad \qquad \qquad \ \ \ \ \ \,\, \textrm{for Euler integrated process},\\
   \rho_W&:=&\liminf_{d\to\infty}\ \frac{r_d}{d^s} \ > \  0 \ \ \ \
   \mbox{for some $s>\tfrac 12$} \  \ \ 
\qquad \textrm{for Wiener integrated process}.
\end{eqnarray*}

For the Korobov case, strong polynomial tractability depends on 
 $\{g_d\}$ and holds iff
$$
  \rho_K:=\liminf_{d\to\infty}\ \frac{\ln\frac1{g_d}}{\ln d}\ >1,
$$

For $g_d=9^{-r_d}$ , we see that $\rho_K=(2\,\ln\,3)\,\rho_E$ and  
conditions for strong polynomial tractability for the Euler and Korobov cases 
are equivalent. 

For $g_d=d^{-r_d/d^s}$, we see that $\rho_W=\rho_K$. Hence,
strong polynomial tractability holds for the Wiener and Korobov 
cases  if $\rho_W>1$, whereas it holds only for the Wiener case when
$\rho_W\in(0,1]$.

Without going to technical details, we may say
that all depends on the two largest eigenvalues 
for the univariate cases. These eigenvalues are quite different
for the Euler and Wiener cases, whereas for the Korobov case they
depend on the weights~$g_k$.  By adjusting these weights,  
the Korobov case  behaves either like the Euler or Wiener case.

\section*{Acknowledgment}

The work of the first and the third authors was done while they  
participated in the Trimester Program
``Analysis and Numerics for
High Dimensional Problems, May-August, 2011, in Bonn, Germany,
and enjoyed warm hospitality of the Hausdorff Research Institute for
Mathematics. 

The work of the first author was supported by
RFBR grants 10-01-00154, 11-01-12104-ofi-m, and by Federal Focused
Program 2010-1.1.-111.128-033.   
The work of the second and third authors was partially supported by
the National Science Foundation.

\vskip 2pc
\noindent{\bf Authors' Addresses:}

\noindent
M. A. Lifshits, 
Department of Mathematics and Mechanics,
St. Petersburg State University,\\
198504 St.Petersburg, Russia, email: mikhail@lifshits.org
\medskip

\noindent
A. Papageorgiou,
Department of Computer Science,
Columbia University,\\
New York, NY 10027, USA, email: ap@cs.columbia.edu
\medskip

\noindent
H. Wo\'zniakowski, Department of Computer Science, Columbia
University,\\
New York, NY 10027, USA, and \\
Institute of Applied Mathematics and Mechanics, University of Warsaw,\\
ul. Banacha 2, 02-097 Warszawa, Poland, email:
henryk@cs.columbia.edu
\end{document}